\documentstyle[eqsecnum,aps,prd]{revtex}
\begin{document}
\draft
\title{On the eigenmodes of compact hyperbolic 3-manifolds}
\author{ Neil J. Cornish \& David N. Spergel}
\address{Department of Astrophysical Sciences, Peyton hall,
Princeton University, Princeton, NJ 08544-1001, USA}

\twocolumn[\hsize\textwidth\columnwidth\hsize\csname
           @twocolumnfalse\endcsname

\maketitle
\widetext
\begin{abstract}
We present a simple algorithm for finding eigenmodes of the Laplacian for
arbitrary compact hyperbolic 3-manifolds. We apply
our algorithm to a sample of twelve manifolds and generate a list
of the lowest eigenvalues. We also display a selection of eigenmodes taken
from the Weeks space.
\end{abstract}
\bigskip
\medskip
]

\narrowtext


Eigenmodes of the Laplace operator contain a wealth of information about
the geometry and topology of a manifold. This is especially true for
hyperbolic 3-manifolds, where the Mostow-Prasad
rigidity theorem\cite{mp} ensures that distinct manifolds have
distinct eigenvalue
spectra: To echo the words of Marc Kac\cite{drum}, one can ``hear the shape''
of a hyperbolic drum.

Unfortunately, the eigenmodes of a compact hyperbolic 3-manifold (CHM)
cannot be expressed in closed analytic form, so numerical solutions must
be sought. A variety of numerical methods exist to solve the problem,
including variational principles based on the finite element method\cite{fe}.
Perhaps the most powerful method is the boundary element method developed
by Aurich and Steiner\cite{as}. Here we present an alternative method
that came out of our work on multi-connected cosmological
models\cite{css1,css2}. While technically inferior to the boundary
element method, our approach is better suited to studying a large
sample of manifolds since the only inputs are the group
generators. In contrast, the boundary element method requires some human
effort prior to each numerical evaluation. To-date, the boundary element
method has only be applied to two 3-dimensional examples, a tetrahedral
orbifold\cite{am}, and the Thurston manifold\cite{ino}. Consequently,
the majority of the eigenvalue spectra described in this paper
are completely new.

\section{Preliminaries}

We seek to solve the eigenvalue problem
\begin{equation}\label{hem}
-\Delta\Psi_q= q^2 \Psi_q\, ,
\end{equation}
for compact hyperbolic 3-manifolds $\Sigma= H^3/\Gamma$, where the
fundamental group, $\Gamma$, is a discrete subgroup of
$SO(3,1)\cong PSL(2,C)$ acting freely and discontinuously.
The metric on the universal cover $H^3$ can be written in
spherical coordinates $(\rho,\,  \theta,\, \phi)$:
\begin{equation}\label{three}
ds^2=d\rho^2+\sinh^2\rho\left(d\theta^2
+\sin^2\theta d\phi^2\right) \, .
\end{equation}
In these coordinates the Laplace operator acting on a scalar function $\Psi$
takes the form
\begin{eqnarray}
&&\Delta \Psi = {1 \over \sinh^2\rho}\left[ {\partial \over \partial \rho}
\left( \sinh^2\rho {\partial \Psi \over \partial \rho} \right) + \right.
\nonumber \\
&&\quad 
\left. {1 \over \sin^2\theta}{\partial \over \partial \theta}\left(
\sin\theta {\partial \Psi \over \partial \theta}\right)
+{1 \over \sin^2\theta}{\partial^2 \Psi \over \partial \theta^2}\right]
\, .
\end{eqnarray}
In the simply connected space $H^3$, the
eigenvalues take all values in the
range $q^2=[1,\infty)$, and the eigenmodes are given by
\begin{equation}
Q_{q\ell m}(\rho,\theta,\phi) = X^{\ell}_{q}(\rho) 
Y_{\ell m}(\theta,\phi) \, .
\end{equation}
Here the $Y_{\ell m}$'s are spherical harmonics and the radial eigenfunctions
are given by the hyperspherical Bessel functions
\begin{equation}
X^{\ell}_{q}(\rho)={ (-1)^{\ell+1} \sinh^\ell\rho\over
 \left(\prod_{n=0}^{l}(n^2+k^2)\right)^{1/2}} 
  \, {d^{\ell+1} \cos(k \rho) \over
d(\cosh\rho)^{\ell+1}} \, .
\end{equation}
The wavenumber, $k$, is related to the eigenvalues of
the Laplacian by
\begin{equation}
k^2=q^2-1 \, .
\end{equation}
The modes have wavelength $\lambda = 2\pi /k$ and an amplitude that
decays as $1/\sinh(\rho)$.
The eigenmodes satisfy the delta-function normalization
\begin{eqnarray}
&&\int_0^\infty \int_0^\pi 
\int_0^{2\pi} \sinh^2\rho d\rho \sin\theta\,  d\theta \, d \phi 
\; \left( Q_{q\ell m} Q^{*}_{q' \ell' m'} \right)\nonumber \\
&& \hspace{0.9in} = \delta(q-q')\delta_{\ell \ell'}\delta_{m m'} \, .
\end{eqnarray}
In principle, the eigenmodes in the multi-connected, compact space
$\Sigma$ can be lifted to the universal cover and expressed in
terms of the eigenmodes of $H^3$:
\begin{equation}
\Psi_q = \sum_{\ell=0}^{\infty}\, \sum_{m=-\ell}^{\ell} a_{q \ell m}
Q_{q\ell m}\, .
\end{equation}
The modes $\Psi_q$ must satisfy the property
\begin{equation}\label{bc}
\Psi_q(x) = \Psi_q(g x) \quad \forall\; g\in \Gamma \; {\rm and}\;
\forall\; x\in H^3 \, ,
\end{equation}
which places restrictions on the expansion coefficients $a_{q\ell m}$.
Indeed, it will only be possible to satisfy (\ref{bc}) when $q^2$ is
an eigenvalue of the compact space. To find the eigenmodes we
numerically solve (\ref{bc}) using a singular value decomposition.

\section{The numerical method}

Our approach to solving (\ref{bc}) is completely straightforward.
We begin by randomly selecting a collection of $d$ points inside the
Dirichlet domain. Taking the face-pairing generators $g_\alpha$
of $\Gamma$, we find all images of our collection of points out to
some distance $\rho_{\rm max}$ in the covering space.
How we chose this distance will be explained later.
Each point $p_j$ yields $n_j$ images and $n_j (n_j+1)/2$ equations
of the form
\begin{equation}
\Psi(g_\alpha p_j) - \Psi(g_\beta p_j) = 0, \quad (\alpha \neq \beta).
\end{equation}
At each point $x=g_\alpha p_j$ the function $\Psi(x)$ is decomposed into
eigenmodes of the covering space with wavenumber $k$:
\begin{equation}
\Psi_k(x) = \sum_{\ell=0}^{L} \sum_{m=-\ell}^{\ell}
a_{k \ell m} Q_{k\ell m}(x)\, .
\end{equation}
How we choose $L$ will be explained in due course. Finding the
$Q_{k \ell m}$'s at each point is easy as there exist numerically
stable recursion relations for both the hyperspherical Bessel
functions and the spherical harmonics.
Schematically, we arrive at the system of equations

\widetext
\begin{equation}\label{main}
\begin{array}{c}
\begin{array}{c}
\uparrow \\
M \\
\downarrow 
\end{array}
\\

\end{array}
\begin{array}{c}
\left[ \begin{array}{ccc}
(Q_{00}(g_1 p_1)-Q_{00}(g_2 p_1)) & \cdots & (Q_{LL}(g_\alpha p_1)
-Q_{LL}(g_\beta p_1)) \\
\vdots  &   \ddots  &  \vdots \\
(Q_{00}(g_1 p_d)-Q_{00}(g_2 p_d)) & \cdots &
(Q_{LL}(g_\alpha p_d)-Q_{LL}(g_\beta p_d))
\end{array}
\right] \\
\longleftarrow \hspace{0.5in} N \hspace{0.5in} \longrightarrow
\end{array}
\begin{array}{c}
\left[
\begin{array}{c}
a_{00} \\
\vdots \\
a_{LL} 
\end{array} \right]
\\

\end{array}
\begin{array}{c}
= \left[ \begin{array}{c}
 0 \\
\vdots \\
0
\end{array}
\right]
\\

\end{array} \; .
\end{equation}
\narrowtext

\noindent The number of rows, $M$, is equal to 
\begin{equation}
M = \sum_{j=1}^d n_j(n_j+1)/2 \, ,
\end{equation}
and the number of columns, $N$, is equal to
\begin{equation}
N = \sum_{\ell=0}^{L} \sum_{m=-L}^{L} 1 \; = (L+1)^2 \, .
\end{equation}
To fix the $a_{klm}$'s up to an overall normalization
requires $M$ to equal $N-1$. When $M=N-1$, a solution exists for any $k$.
However, if $M > N$
a solution only exists when $k$ corresponds to an eigenvalue of
the compact space. The standard numerical method for solving
over-constrained systems of equations is the Singular Value
Decomposition (SVD)\cite{nreps}. For a system of the form
\begin{equation}
\overline{\overline{A}}\cdot \overline{a} = \overline{0} \, ,
\end{equation}
the SVD returns solution vectors $\overline{a}$ that minimize
\begin{equation}
\chi^2 = \left| \,  \overline{\overline{A}}\cdot
\overline{a} \, \right|^2 \, .
\end{equation}
By incrementing in $k$, the eigenvalues are revealed by minima in
the function $\chi^2(k)$. Eigenmodes with multiplicity greater than
unity will yield multiple solution vectors $\overline{a}$.

All that remains to be done is to decide on optimal choices for
$L=\ell_{\rm max}$, the tiling radius $\rho_{\rm max}$ and the degree of
over-constraint $c=M/N$. The choice of $L$ and $\rho_{\rm max}$ is
dictated by the structure of the radial eigenfunctions
$X^{\ell}_{k}(\rho)$. In broad outline, they are of the form
\begin{equation}
X^{\ell}_{k}(\rho) \approx \left\{ \begin{array}{ll}
0 & \rho \ll \rho_0 \\
\\
\displaystyle{ \frac{ \cos(k\rho +\phi_0)}{ \sinh(\rho)}} &
\rho \gg \rho_0 \, .
\end{array} \right.
\end{equation}
The constants $\phi_0$ and $\rho_0$ depend on $k$ and $\ell$. For
fixed $k$, $\rho_0$ increases monotonically with increasing $\ell$.
Therefore, if we restrict our attention to some finite region with
$\rho \leq \rho_0$, we need only consider a finite number
of multipoles $L$. With $\ell$ held fixed, $\rho_0$ decreases
monotonically with increasing $k$. Thus, if we hold
$\rho_{\rm max}\approx \rho_0$
fixed, we must increase $L$ as $k$ increases. Alternatively, if we
hold $L$ fixed, we must decrease $\rho_{\rm max}$ as $k$ increases.

Because the number of computational steps scales as $L^6$, it makes
sense to keep $L$ as small as possible. However, small values of
$L$ yield small values of $\rho_0$, which

\medskip

\hrule

\vspace{1.2in}

\hrule

\medskip

\noindent in turn limits the number of
images we are able to collect for each point. As a compromise, we
choose $L$ to be as small as possible, while leaving $\rho_{\rm max}$
large enough for there to be at least 10 images of each point.
Numerically we set $\rho_{\rm max}$ to be the first solution to
the transcendental equation
\begin{equation}
X^{L}_{k}(\rho)\, \sinh(\rho) = 0.25 \, .
\end{equation}
We also found it advantageous to fix a minimum radius in the same way,
but with $L$ replaced by $\ell_{\rm min}$ in the above equation. The
inner cut-off helps to keep the $Q_{k\ell m}$'s of similar size.
Finally, we found the optimal degree of oversampling to lie in the
range $10 \rightarrow 100$, with the high end of the range being required
at low $k$ and $L$. At higher $k$ the radial eigenfunctions look less
alike and a lower oversampling can be used. Good all purpose choices
covering the range $k=1$ to $k=20$ are:
\begin{eqnarray}
&& L = 10 + [k] \, , \nonumber \\
&& \ell_{\rm min} = 5 \, , \nonumber \\
&& c = 10 + [100/k] \, ,
\end{eqnarray}
though the algorithm performs well for a wide range of inputs.
The eigenmodes derived using the above choice of inputs are good
to within a few percent. Choosing $L$ and $c$ larger improves the
eigenmodes, but at the cost of slowing down the computations.

The SVD returns solution vectors for the $a_{k\ell m}$ normalized such
that
\begin{equation}
\sum_{\ell=0}^{L} \sum_{m=-\ell}^{\ell}
\left| a_{k \ell m} \right|^2 = 1 \, .
\end{equation}
Thus, the eigenmodes are automatically delta-function normalized in $H^3$.
To normalize the modes in the compact space $\Sigma$ we numerically
perform the integral
\begin{equation}\label{norm}
\int_{\Sigma} \Psi_k(x) \Psi^*_{k'}(x)\, dV \, .
\end{equation}
When $k=k'$ we normalize to one, and when $k\neq k'$ we check that
the integral vanishes (or is at least tolerably small).

\subsection{Example}

We begin by choosing an example from Jeff Weeks' {\em SnapPea}\cite{weeks}
census of closed hyperbolic 3-manifolds and ask {\em SnapPea} for
the face-pairing generators. Taking $m188(-1,1)$ as a randomly chosen
example, we ran our code out to $k=10$ to produce the $\chi^2(k)$ shown
in Fig.~1. The $\chi^2$'s for the five best solution vectors are shown.
The eigenvalues appear as clear minima in the $\chi^2$ of the
best fit solution vector. The first mode to have multiplicity
greater than one occurs at $k=8.34$. However, below $k=8$ we see that
the $\chi^2$ of the second best solution vector occasionally
exhibits minima at local maxima of the best fit solution vector. This occurs
when the second best fit is formed from a weighted
superposition of adjacent eigenmodes. In contrast, when an eigenmode has
multiplicity greater than one the minima of the second best fit coincides
with that of the best fit. This behaviour is due to the SVD returning an
orthonormal set of solution vectors. 

When successive minima are close
together, such as occurs near $k=9.8$ in Fig.~1, care must be taken in
locating successive minima as they can be displaced from their true
positions. This problem is a familiar one for astronomers who
study spectral lines in starlight, and we were able to use the same
deconvolution techniques to study our eigenspectra. The results of
our analysis are compiled in Table I, where the eigenvalues and their
multiplicities are recorded.

\begin{table}
\caption{Eigenvalue spectrum, $q^2$, for m188(-1,1).}
\begin{tabular}{cccccccc}
\\
& 20.4 & 22.6 & 27.2 & 30.2 & 39.6 & 46.2 & \\
& 1    &  1   &  1   &  1   &  1   &  1   & \\
\hline
\\
& 51.8 & 55.3 & 60.1 & 70.6 & 75.5 & 78.8 & \\
&  1   &   1  &  1   &  2   &   2  &  1   &  \\
\hline
\\
& 80.9 & 83.1 & 86.0 & 96.8 & 98.0 & 99.4 & \\
&  1   &   1  &  1   &  2   &   1  &  1   & \\
\end{tabular}
\end{table}

\
\begin{figure}[h]
\vspace{50mm}
\includegraphics{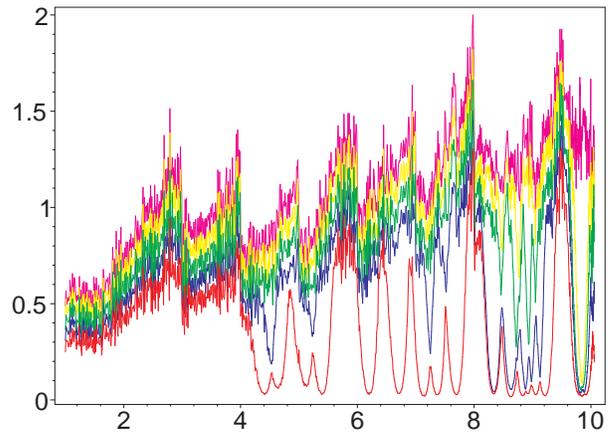}
\vspace{1mm}
\caption{The $\chi^2$ spectrum for $m188(-1,1)$ in the range
$k=1 \rightarrow 10$.}
\end{figure}

\section{Results and checks}

We applied a battery of tests to our results. The first
was to reproduce Inoue's\cite{ino} results for the Thurston manifold.
Calling up $m003(-2,3)$ from the {\em SnapPea} census, we generated the
collection of low-lying eigenmodes recorded in Table II.

Inoue was able to generate modes with $q^2 \leq 100$, and in this range
our eigenvalues agree. However, we discovered that Inoue had missed
eigenmodes at $q^2=46.2$ and $q^2=59.1$.
Past $q^2=100$ we are in uncharted territory, so other checks have to be
applied. One simple check is to compare the spectral staircase, {\it i.e.}
the number of modes with wavenumber $\leq k$, with the prediction
of Weyl's asymptotic formula
\begin{equation}
N(\leq k) \asymp \displaystyle{ \frac{{\rm Vol}(\Sigma)}{6 \pi^2}}k^3
+ {\rm const}. \; .
\end{equation}

\begin{table}
\caption{Eigenvalue spectrum for the Thurston space.}
\begin{tabular}{cccccccc}
 \\
& 29.3 & 33.5 & 46.2 & 47.8 & 50.8 & 59.1 & \\
& 1   &   1   &   2  &  1   &  1   &  2  & \\
\hline
\\
& 68.9 & 73.8 & 76.2 & 85.8 & 95.1 & 98.0 & \\
&  1   &   1  &  1   &  2   &  1   &  1   & \\
\hline
\\
& 100.1 & 107.5 & 113.8 & 115.5 & 117.4 & 123.3 & \\
&  1   &  1     &   1   &   2   &   1   &  2    &  \\
\hline
\\
& 130.3 & 137.9 & 140.0 & 144.5 & 149.8 & 156.3 & \\
&   1   &   2   &   1   &   2   &   1   &   1   & \\
\hline
\\
& 160.5 & 164.1 & 166.6 & 169.5 & 175.0 & 178.2 & \\
&  1    &   2   &   2   &   1   &   1   &   1   &  \\
\hline
\\
& 184.6 & 192.3 & 197.8 & 204.3 & 207.2 & 209.8 \\
&   2   &   2   &   4   &  2    &   1   &   1   \\
\end{tabular}
\end{table}

From Fig.~2 we see that the spectral staircase closely follows the Weyl
formula. Another check we applied to the eigenmodes
was to evaluate the integral (\ref{norm}) for modes with
unequal $k$. The typical overlap was found to be on the order
of a few percent, which is comparable to the performance of the
boundary element method. 

\
\begin{figure}[h]
\vspace{60mm}
\includegraphics{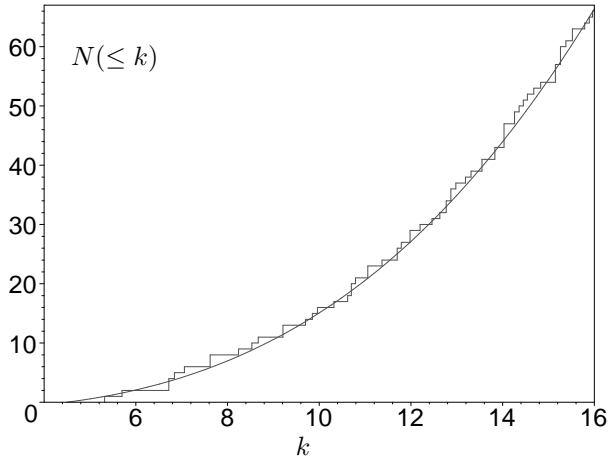}
\vspace{1mm}
\caption{The spectra staircase for $m003(-2,3)$ showing its convergence toward
the Weyl asymptotic formula.}
\end{figure}
\vspace*{-2mm}
\begin{picture}(0,0)
\put(30,190){$N(\leq k)$}
\put(115,42){$k$}
\end{picture}

\subsection{The Weeks Space}

To really put our code to the test we decided to study the Weeks space.
This space has the distinction of being the smallest known hyperbolic
3-manifold, and is now generally thought to be {\em the} smallest
example. The Weeks space provides a challenge to eigenmode solvers as it
has an unusually large symmetry group. Consequently, many of the modes
will be highly degenerate. The symmetry group is the Dihedral group of
order 6, which has the presentation
\begin{equation}
D_6 = < a, b \; \vert \, a^6, \, b^2, \, ba^{-1}ba^{-1} > \, .
\end{equation}
The geometrical interpretation is that there exists a closed geodesic
about which the space has a six-fold rotational symmetry, and a second
closed geodesic, orthogonal to the first, about which the manifold
has a reflection symmetry. If we were to choose the basepoint of our
Dirichlet domain at one of the points where the two geodesics intersect,
the resulting fundamental polyhedron would enjoy the full $D_6$
symmetry\footnote{We thank Jeff Weeks and Craig Hodgson for providing
us with this description.}.
Because of the $D_6$ symmetry, we expect to find modes with 1, 2, 3, 4 and
6 fold degeneracy.

Taking the Weeks space $m003(-3,1)$ from the {\em SnapPea} census, we
located the first 74 eigenmodes. These are listed in Table III
along with their multiplicities. Many of the higher eigenmodes are
highly degenerate, in keeping with our expectations. The spectral staircase
shown in Fig.~3 is in excellent agreement with Weyl's asymptotic formula.

For those curious to see what the eigenmodes themselves look like, we
display a series of slices through a selection of modes. These appear in
Figures 5 and 6. For reference we also display a view of the Dirichlet
domain in Fig.~4, to help make contact with the 3-dimensional structure
of the modes.

\
\begin{figure}[h]
\vspace{60mm}
\includegraphics{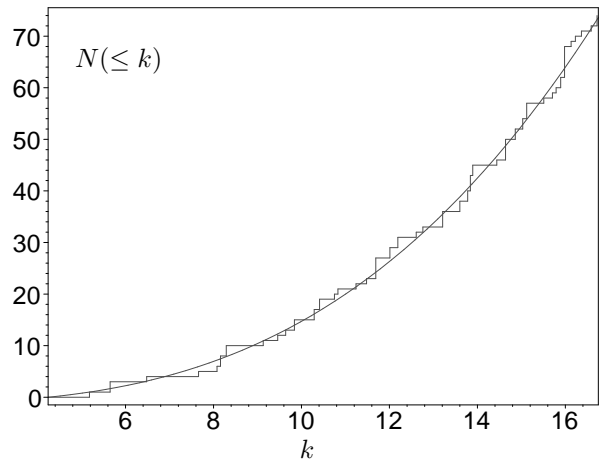}
\vspace{1mm}
\caption{The spectra staircase for $m003(-3,1)$ showing its convergence toward
the Weyl asymptotic formula.}
\end{figure}
\vspace*{-2mm}
\begin{picture}(0,0)
\put(30,190){$N(\leq k)$}
\put(115,42){$k$}
\end{picture}

\begin{table}
\caption{Eigenvalue spectrum for the Weeks space.}
\begin{tabular}{cccccccc}
 \\
& 27.8 & 32.9 & 43.0 & 59.7 & 66.3 & 67.6 & \\
& 1   &   2  &   1  &  1   &  1   &  2   & \\
\hline
\\
& 69.7 & 84.4 & 90.5 & 93.9 & 97.8 &  106.9 & \\
& 2    &   1  &  1   &  1   &   2  &   2  & \\
\hline
\\
& 109.4 & 116.6 & 118.3 & 127.3 & 132.8 & 137.7 & \\
&  2    &   1   &   1   &   1   &   1   &   4   & \\
\hline
\\
& 145.2 & 149.6 & 160.0 & 163.8 & 175.5 & 186.0  & \\
&   2   &   2   &   1   &   1   &   3   &   2    & \\
\hline
\\
& 190.9 & 192.6 & 194.1 & 209.5 & 215.3 & 221.8 & \\
&  2    &   3   &  2    &  1    &   4   &   2   & \\
\hline
\\
& 226.9 & 229.6 & 241.6 & 247.6 & 250.3 & 253.6 & \\
&  2    &   3   &   1   &   1   &   1   &   2   &  \\
\hline
\\
& 256.4 & 261.2 & 264.3 & 268.8 & 276.2  & 280.6 & \\
&  6    &   1   &   1   &   1   &   1    &   2   &  \\
\end{tabular}
\end{table}

\newpage

\
\begin{figure}[h]
\vspace{70mm}
\includegraphics{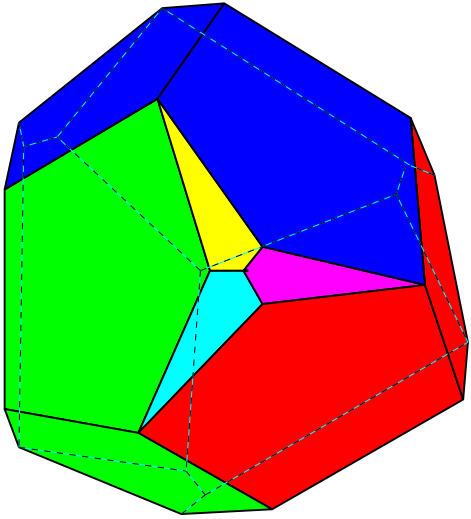}
\vspace{1mm}
\caption{A Dirichlet domain for the Weeks space shown in Klein coordinates.
Like-colored sides of the polyhedron are topologically identified.}
\end{figure}

\begin{figure}[h]
\vspace{100mm}
\includegraphics{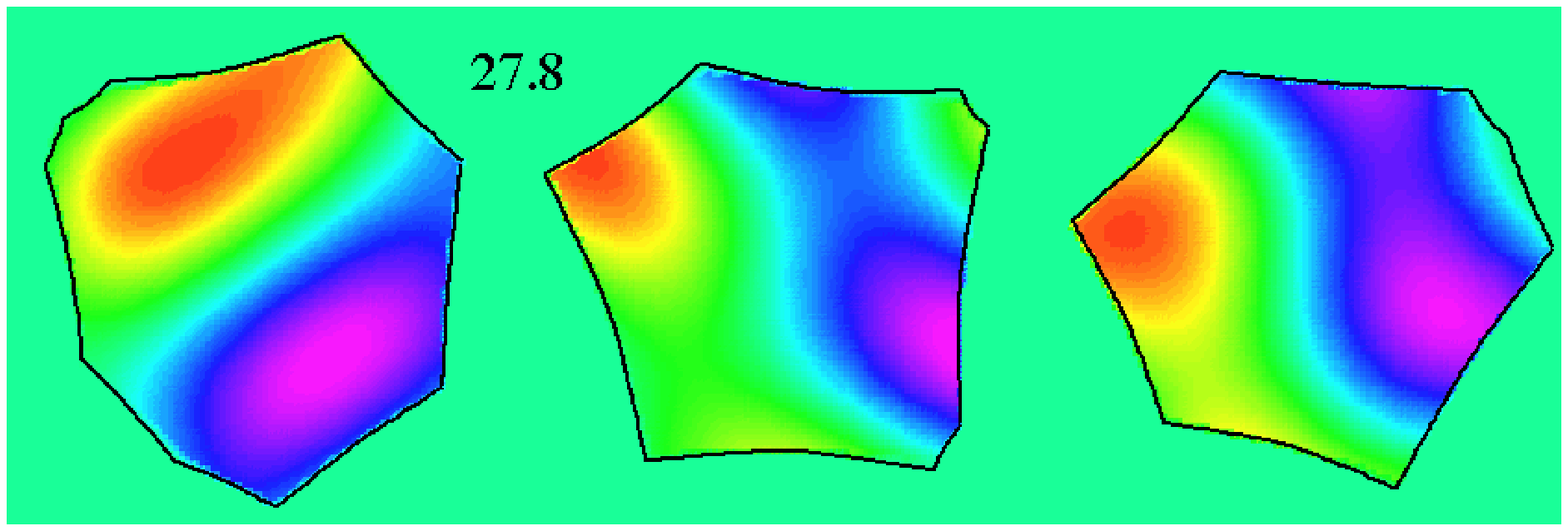}
\includegraphics{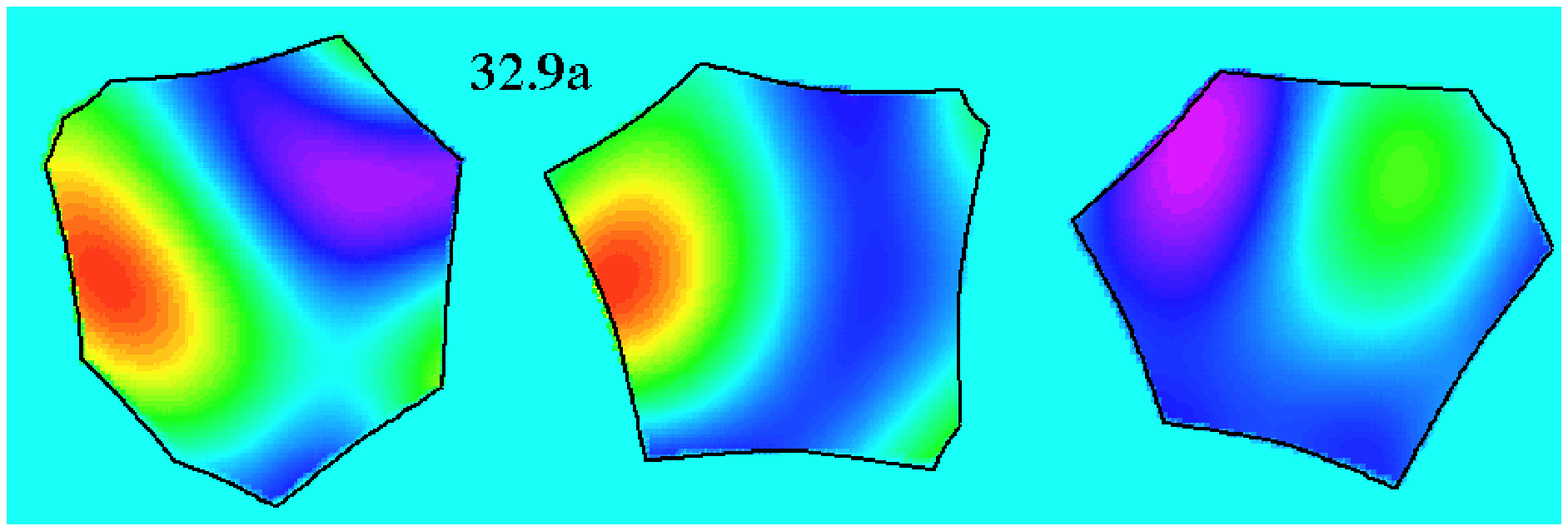}
\includegraphics{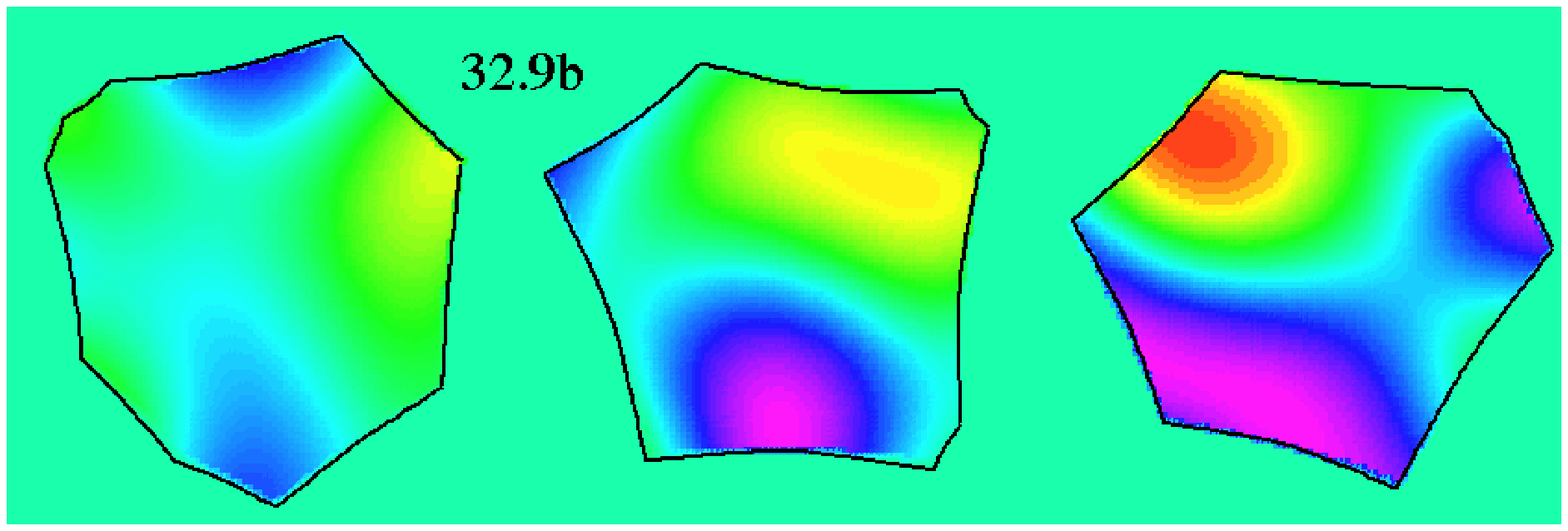}
\vspace{1mm}
\caption{The first three eigenmodes of the Weeks space. The three
views in each panel are, respectively, the $x=0$, $y=0$ and $z=0$ slices
through the fundamental domain. Here we are using Poincar\'{e} coordinates.}
\end{figure}

\newpage

\
\begin{figure}[h]
\vspace{95mm}
\includegraphics{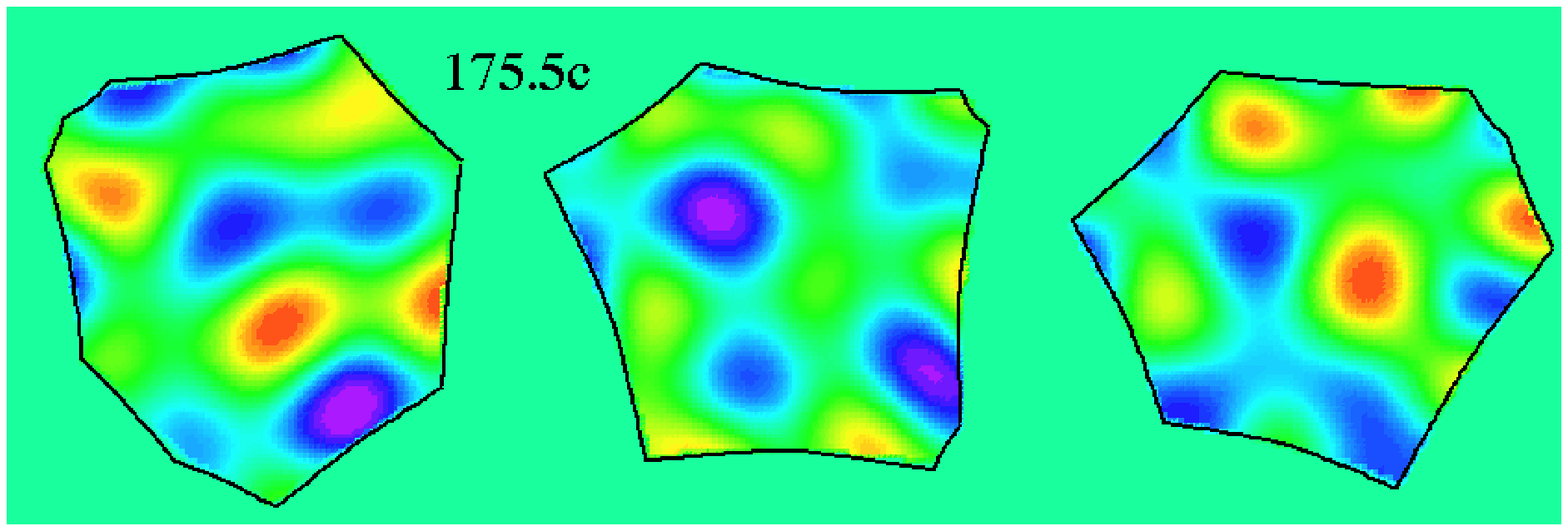}
\includegraphics{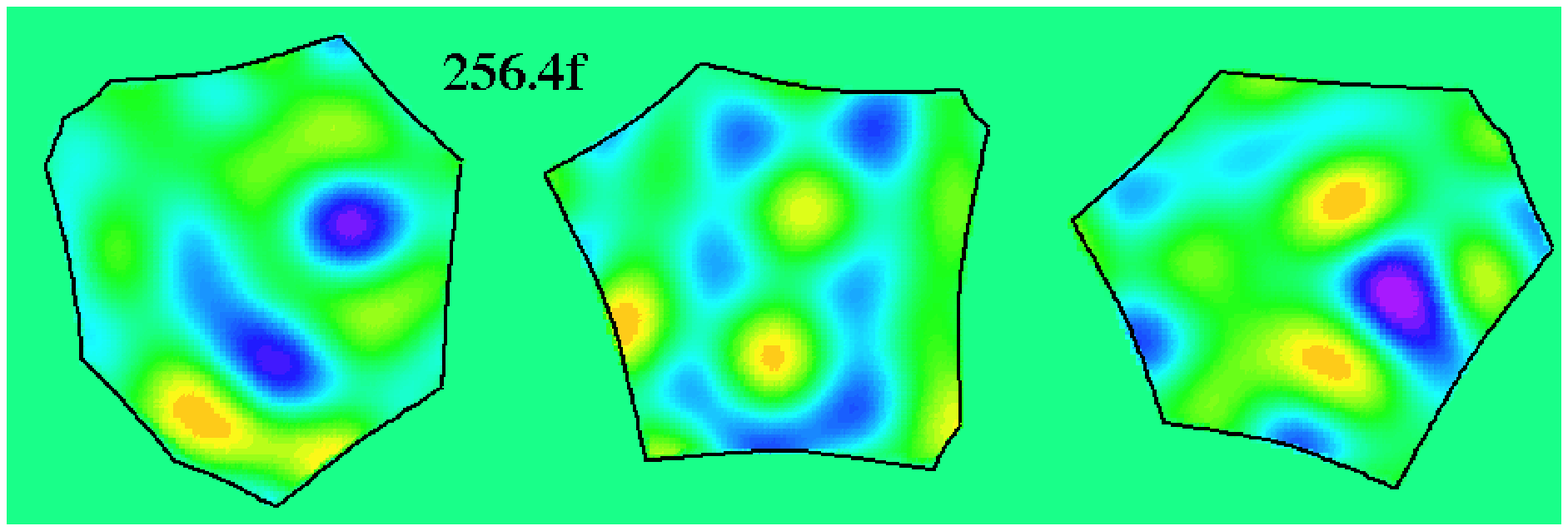}
\includegraphics{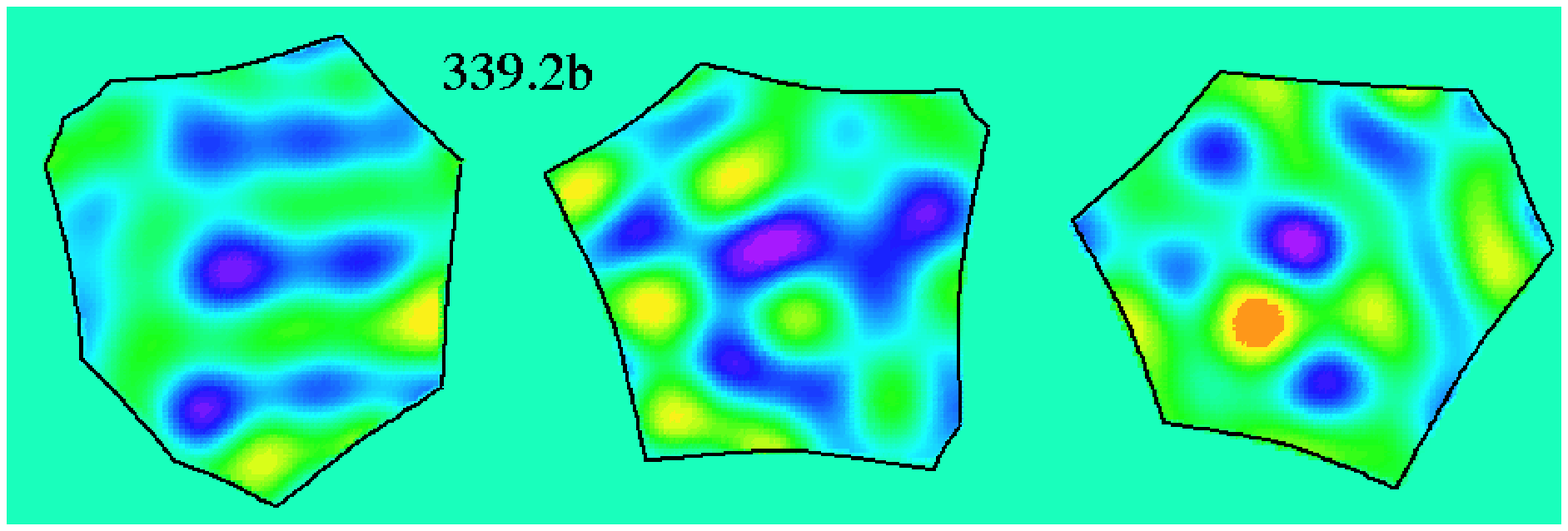}
\vspace{1mm}
\caption{The first of the 3 and 6 fold degenerate modes
of the Weeks space. Also shown is the highest mode we generated.}
\end{figure}

\subsection{The GOE Prediction}

Compact hyperbolic manifolds provide the archetypal setting for
chaos. Consequently, we expect the statistical properties of the modes to
be described by random matrix theory\cite{rmt}. Because the modes are
associated with time-reversal invariant dynamics, we expect
the statistical properties to be those of a Gaussian Orthogonal
Ensemble (GOE). The GOE prediction is that the quantity
\begin{equation}
x = \displaystyle{ \frac{ |a_{k\ell m}-\bar{a}_k|^2}{\sigma_k^2}} \, ,
\end{equation}
should behave as a pseudo-random number with probability distribution
\begin{equation}
P(x) = \displaystyle{\frac{1}{\sqrt{2 \pi x}}}\, e^{-x/2} \, .
\end{equation}
In the above equations, $\bar{a}_k$ denotes the average of the $a_{k\ell m}$'s
and $\sigma^2_k$ their variance. To avoid the singularity at $x=0$, it
is conventional to compare numerical results to the cumulative
distribution
\begin{equation}
I(x) = \int_0^x P(x)\, dx = {\rm erf}(\sqrt{x/2}) \, .
\end{equation}

Taking the mode with eigenvalue $q^2=175.5$ as an example, we display
in Fig.~7 the first 676 $a_{k\ell m}$'s in a scatter plot. Notice that the
distribution is independent of $\ell$ and $m$, as expected for a
chaotically mixed state.

\
\begin{figure}[h]
\vspace{60mm}
\includegraphics{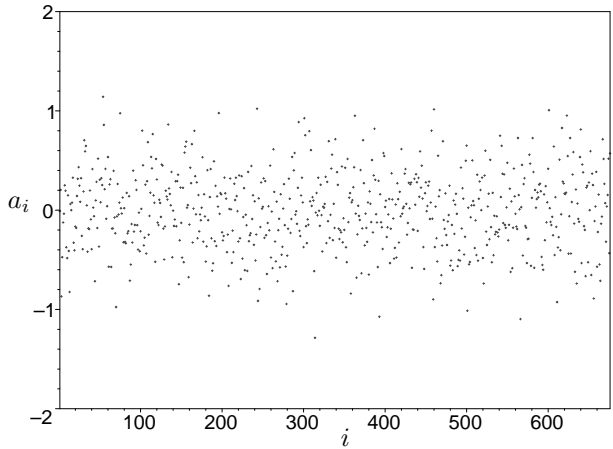}
\vspace{1mm}
\caption{The $a_{k\ell m}$'s for the eigenmode $q^2=175.5$, shown up to
$\ell_{\rm max} = 25$. The numbering scheme for the modes maps
$a_{\ell m} \mapsto a_{i}$ where $i=\ell^2+\ell+1+m$.}
\end{figure}
\vspace*{-2mm}
\begin{picture}(0,0)
\put(-5,150){$a_{i}$}
\put(121,60){$i$}
\end{picture}

The cumulative distribution $I(x)$ is compared to the GOE prediction
in Fig.~8. The agreement is quite remarkable.

\
\begin{figure}[h]
\vspace{60mm}
\includegraphics{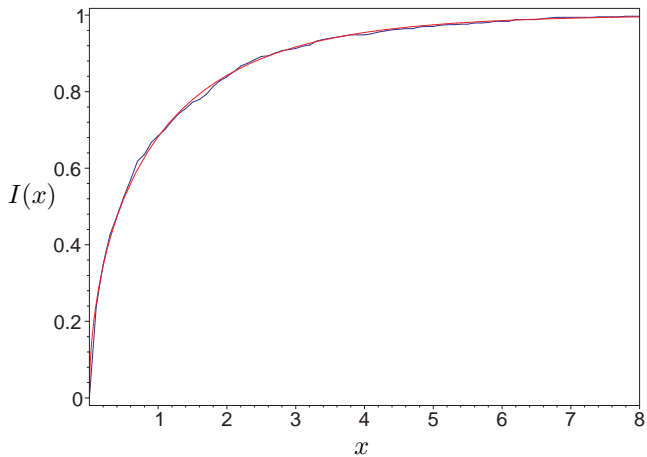}
\vspace{1mm}
\caption{The cumulative distribution $I(x)$ found for the eigenmode
with $q^2=175.5$ (blue line) and the GOE prediction (red line).}
\end{figure}
\vspace*{-2mm}
\begin{picture}(0,0)
\put(-6,150){$I(x)$}
\put(125,55){$x$}
\end{picture}

\subsection{The circles test}

The next test we applied to the modes is one that we plan to
apply to our own universe\cite{css1,css2,jw}. Imagine drawing a
2-sphere of radius $\rho$ about the basepoint of the Dirichlet
domain. Viewed in the universal cover, the space will contain
an infinite number of copies of this 2-sphere. If the radius of
the 2-sphere exceeds the in-radius of the Dirichlet domain, then
the 2-spheres will intersect along a circle. Mapping the entire
picture back inside the Dirichlet domain, we see that the 2-sphere
self-intersects. If we now take a 2-sphere slice through one of the
eigenmodes, the amplitude of the mode must match up around the matched
pair of circles.

The largest matched circles lie on face-planes of the Dirichlet
domain. Taking a face-pairing generator $g$, and representing it
as a $4\times 4$ real matrix in $SO(3,1)$, the angular radius $\alpha$
of the matched circle is
\begin{equation}
\alpha = {\rm arccos}\left( \displaystyle{\frac{g_{00}-1}
{\sqrt{g_{00}^2-1}\tanh(\rho)}} \right) \, .
\end{equation}
The $(\theta,\phi)$ coordinates of the circle centers are
\begin{eqnarray}
&& \left({\rm arccos}\left(g_{10}/\sqrt{g_{00}^2-1}\right),
{\rm arctan}\left(g_{30}/g_{20}\right)\right)\, , \nonumber \\
&& \left({\rm arccos}\left(g^{-1}_{10}/\sqrt{g_{00}^2-1}\right),
{\rm arctan}\left(g^{-1}_{30}/g^{-1}_{20}\right)\right) \; .
\end{eqnarray}

We display in Fig.~9 the amplitude of the eigenmodes $q^2=43.0$ and
$q^2=84.4$ around a 2-sphere of radius $\rho=1$
using an equal-area projection.
Two of the matched circle pairs are indicated by white lines.
In Fig.~10 we plot the amplitude of each mode around each pair
of circles and see that they are indeed properly matched. If our
universe is multi-connected, we hope to see similar matched circles in
the cosmic microwave sky\cite{css1,css2}.

\
\begin{figure}[h]
\vspace{80mm}
\includegraphics{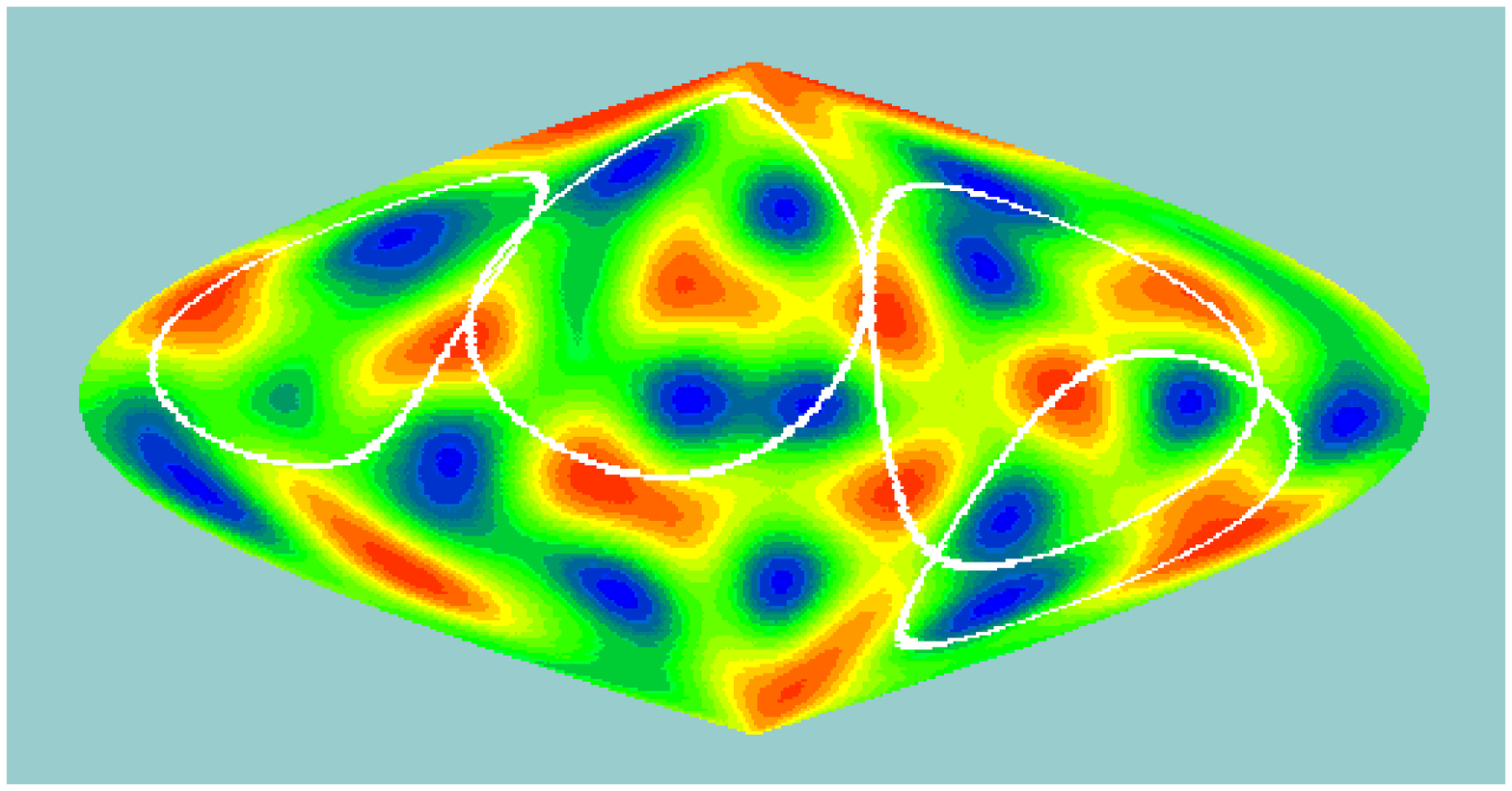}
\includegraphics{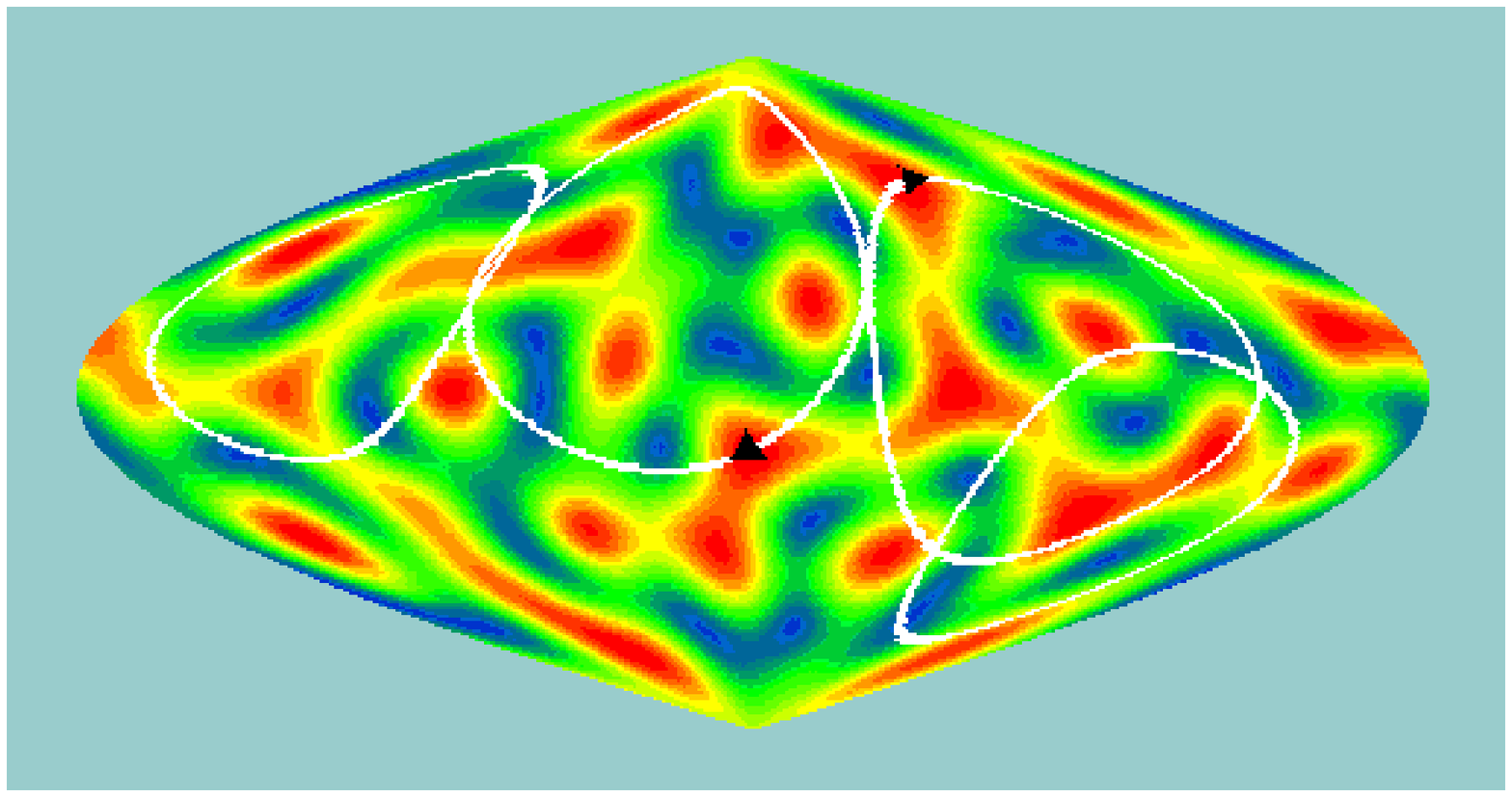}
\vspace{8mm}
\caption{A two-dimensional slice through the $q^2=43.0$ (upper panel)
and $q^2=84.4$ (lower panel) eigenmodes. The slice is taken on a 2-sphere
of unit radius. Two pairs of matched circles are also indicated. The
relative phasing of one matched pair is shown on the lower panel.}
\end{figure}

\newpage

\
\begin{figure}[h]
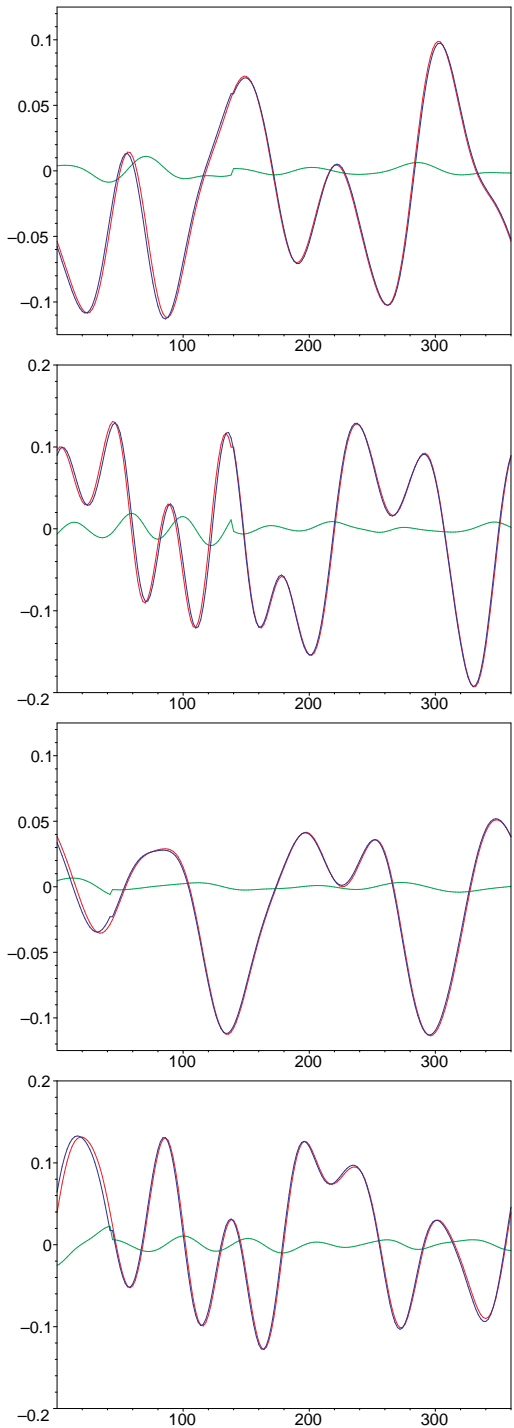

\vspace{187mm}
\includegraphics{weeks648_p5.ps}
\includegraphics{weeks913_p5.ps}
\includegraphics{weeks648_p8.ps}
\includegraphics{weeks913_p8.ps}
\vspace{4mm}
\caption{ Each panel shows $\Psi_q$ around a pair of matched circles, $C$
(in red) and $\widetilde{C}$ (in blue). Also shown (in green)
is the difference $\Psi_q-\widetilde{\Psi}_q$
around the matched pair. This small difference represent the
numerical error in our algorithm. The first panel correspond to
the mode $q^2=43.0$ sampled along the inner pair of circles. The second panel
shows the same pair of circles, but for the mode $q^2=84.4$. The third and
fourth panels are likewise shown for
the outer pair of circles. }
\end{figure}

\section{The Lowest Eigenmodes}

Having established that our algorithm is reliable, we set it loose
on the {\em SnapPea} census to produce the list of lowest eigenvalues
recorded in Table IV. The lowest eigenvalue, $q^2_1$, is a useful topological
invariant that has attracted considerable attention in the mathematical
literature. There exist a variety of upper and lower bounds on $q_1^2$.
A summary of these bounds can be found in the works of Callahan\cite{pat}
and Cornish {\it et al.}\cite{css3}. Amongst the sharpest are those that
employ the diameter\footnote{The diameter is defined to equal the
greatest distance between any two points in the manifold.}, $D$, of the space:
\begin{equation}\label{bounds}
\displaystyle{ \frac{ 4\widetilde{D}}{ D^2 (\sinh(\widetilde{D})+
\widetilde{D})^2}} \leq \; q_1^2 \; \leq 1+\left(\displaystyle{\frac{2\pi}{D}}
\right)^2 \, .
\end{equation}
Here $\widetilde{D}$ denotes the square root of the smallest integer
that is greater than or equal to $D^2$. A listing of the diameters can
be found in Table V. In all cases, the lowest eigenvalues found by our
algorithm fell in the range dictated by (\ref{bounds}). The eigenvalue
bound (\ref{bounds}) tells us that the wavelength of the lowest
eigenmode must be greater than or equal to the diameter of the space.
Curiously, we found that $\lambda_1 = (1.3 \rightarrow 1.6) D$ for all twelve
examples studied.

On a cautionary note, the eigenvalues listed in Table IV might not
be the lowest supported by these spaces. Our method for finding the
eigenmodes is unable to detect modes with $q^2 < 1$, as these
modes have imaginary wavenumbers. Moreover, the spherical eigenmodes
we use as our expansion basis start to look very much alike for
wavenumbers $k<1$, so we may have missed modes in the range
$k=[0,1/4]$. We are currently developing a variant of our method to
handle all modes below $q^2 = 2$. At the other end of the spectrum,
the only limitation in going out to $q^2=\infty$ is computer power.
To get all the 

\begin{table}
\caption{Lowest eigenvalues}
\begin{tabular}{ccccc}
$\Sigma$& Vol& $G$ & $q^2_1$ & $m_1$ \\
\hline 
m003(-3,1) & 0.9427 & $D_6$ & 27.8 & 1 \\
m003(-2,3) & 0.9814 & $D_2$ & 29.3 & 1 \\
s556(-1,1) & 1.0156 & $Z_4$ & 27.9 & 1 \\
m006(-1,2) & 1.2637 & $D_4$ & 21.1 & 2 \\
m188(-1,1) & 1.2845 & $D_2$ & 20.4 & 1 \\
v2030(1,1) & 1.3956 & $D_2$ & 16.2 & 1 \\
m015(4,1)  & 1.4124 & $D_2$ & 28.1 & 2 \\
s718(1,1)  & 2.2726 & $D_2$ & 10.1 & 1 \\
m120(-6,1) & 3.1411 & $Z_2$ & 7.50 & 1 \\
s654(-3,1) & 4.0855 & $D_2$ & 5.88 & 1 \\
v2833(2,3) & 5.0629 & $Z_2$ & 6.29 & 1 \\
v3509(4,3) & 6.2392 & $D_2$ & 6.06 & 1
\end{tabular}
\end{table}

\begin{table}
\caption{$L_\gamma$, diameter, wavelength and wavenumber}
\begin{tabular}{ccccc}
$\Sigma$& $L_\gamma$ & $D$ & $\lambda_1/D$ & $k_1$ \\
\hline 
m003(-3,1) & 0.585 & 0.843 & 1.44  & 5.18 \\
m003(-2,3) & 0.578 & 0.868 & 1.36  & 5.32 \\
s556(-1,1) & 0.831 & 0.833 & 1.45  & 5.19 \\
m006(-1,2) & 0.575 & 1.017 & 1.38  & 4.48 \\
m188(-1,1) & 0.480 & 0.995 & 1.44  & 4.41 \\
v2030(1,1) & 0.366 & 1.082 & 1.49  & 3.90 \\
m015(4,1)  & 0.794 & 0.923 & 1.31  & 5.21 \\
s718(1,1)  & 0.339 & 1.439 & 1.45  & 3.01 \\
m120(-6,1) & 0.314 & 1.694 & 1.45  & 2.55 \\
s654(-3,1) & 0.312 & 1.946 & 1.46  & 2.21 \\
v2833(2,3) & 0.486 & 1.701 & 1.60  & 2.30 \\
v3509(4,3) & 0.346 & 1.802 & 1.55  & 2.25
\end{tabular}
\end{table}

\noindent modes out to $q^2=250$ takes around 12 hours on a
single R1000 Silicon Graphics chip (or 10 minutes if you use
a 64 processor Origin 2000 supercomputer as we did).

\section{Acknowledgments}
We are indebted to Ralf Aurich, Craig Hodgson and Jeff Weeks for
answering our many questions about the structure of the eigenmodes,
the properties of the symmetry groups and the topology of 3-manifolds.
We are grateful for the support provided by NASA through their funding
of the MAP satellite mission {\tt http://map.gsfc.nasa.gov/}.

\end{document}